\documentclass[12pt]{amsart}
\usepackage[margin=1.0in]{geometry}
\usepackage{mathptmx,color}
\usepackage[foot]{amsaddr}
\begin{document}
\title{Geometry in the Courtroom}
\author{Noah Giansiracusa$^*$} 
\address{$^*$Assistant Professor of Mathematics, Swarthmore College} 
\author{Cameron Ricciardi$^{**}$} 
\address{$^{**}$Undergraduate Math Major, Swarthmore College}
\email{ngiansi1@swarthmore.edu, criccia2@swarthmore.edu}

\maketitle

\begin{abstract}
There has been a recent media blitz on a cohort of mathematicians valiantly working to fix America's democratic system by combatting gerrymandering with geometry.   While statistics commonly features in the courtroom (forensics, DNA analysis, etc.), the gerrymandering news raises a natural question: in what \emph{other} ways has pure math, specifically geometry and topology, been involved in court cases and legal scholarship?  In this survey article, we collect a few examples with topics ranging from the Pythagorean formula to the Ham Sandwich Theorem, and we discuss some jurists' perspectives on geometric reasoning in the legal realm.  One of our goals is to provide math educators with engaging real-world instances of some abstract geometric concepts.
\end{abstract}

\section{Introduction}

With the upcoming national census in 2020 and ensuing congressional redistricting, the problem of gerrymandering has become a tremendously hot issue, and rightly so.  One of the legal guidelines courts have established is that districts must be ``compact,'' though no proper definition of this term has been provided (all districts are closed and bounded, so the one you're thinking of won't help!).  Occasionally a particularly egregious district is struck down based on ad hoc arguments; due to the scope and complexity of the problem, however, a more scientific approach is clearly required.  Toward this end, mathematician Moon Duchin formed a group\footnote{See \texttt{https://sites.tufts.edu/gerrymandr/}} embarking on an ambitious interdisciplinary project with lawyers, political scientists, geospatial analysts, and computer scientists to apply rigorous, sophisticated geometric methods to the problem of detecting and preventing gerrymandering.  Press coverage of these efforts has appeared in \emph{Nature}, \emph{Wired}, \emph{NBC News}, \emph{NPR}, etc., often with superficial headlines implicitly suggesting bewilderment that esoteric mathematics may actually be useful in the real world (gasp!).  

While witnessing this deluge of public attention, a simple question came to our minds: when else has geometry entered the courtrooms, either in actual cases or in legal scholarship?  After finding some intriguing and entertaining examples, we decided to write this survey paper to collect these in one place and share them with the mathematics community.  We hope math teachers will find ways to incorporate some of this material into their curricula.\\

\emph{Acknowledgements}.  We would like to thank Professor Paul Edelman\footnote{Professor of Mathematics and Professor of Law at Vanderbilt University.} for encouragement and help finding references, and Professor Talvacchia\footnote{Professor of Mathematics at Swarthmore College.} for helpful conversations pertaining to \S\ref{sec:maps}.  The first author was supported in part by an NSA Young Investigator Grant and wishes to thank that organization for financial support.

\section{Background}

There was a famous debate in the \emph{Harvard Law Review} in the early 1970s on what role, if any, mathematics and quantitative methods should play in jurisprudence \cite{Finkelstein-Fairley,Tribe,Finkelstein-Fairley2,Tribe2} (see also \cite{Kaye,MathTrial,ConCal,Noah-law}).  This remains a difficult topic, and one we shall not address here.  Despite caution and controversy, innovation in the technological and scientific realms has forced the courtroom to accept certain forms of mathematically-framed evidence (think forensics and DNA analysis, among many other topics \cite{StatsLawyers,MathForensics}), which generally fall under the purview of statistics. Given the modern-day ubiquity of statistics, this really comes as no surprise,\footnote{This is not to say that statistical evidence in the courtroom is a uniquely modern development; indeed, there is a fascinating case from 1868 in which Harvard Math Professor Benjamin Peirce used the binomial distribution to argue that the signature on a disputed amendment to a wealthy matron's will was a forgery \cite{StatsWill}.} but what was (for us, at least) less expected is that various geometric concepts from supposedly ``pure'' mathematics would find their way into the legal realm, and it is this that we focus on in this paper.\footnote{Since there is presently a surge of interest and writing in the geometry of gerrymandering \cite{GerryNature,GerryAMS,GerryGap,GerryCompactness,GerryConvexity}, we mostly avoid this topic in order to give attention to other instances of geometry in the courtroom.}

In some situations geometric considerations directly enter the legal arguments and judges' rulings, while in other situations a geometric analysis is provided by an academic outside the legal case.  Moreover, some in the legal profession have used geometric reasoning abstractly to better understand delicate legal issues, whereas others have argued against the idea of viewing law as an axiomatic discipline modeled on Euclidean geometry.  One of the main matters that arises here is that of precision.  Mathematics is founded on precisely stating all axioms, definition, theorems, etc., and avoiding ambiguity like the plague.  Law, on the other hand, necessarily deals with amorphous concepts that defy precise characterization; in fact, it is often advantageous to remain deliberately vague so as to allow a natural development (cf., \cite{Vagueness}).  For instance, if the courts had rigorously defined what they meant by ``compact'' in the setting of districting, then the prospect is raised that a duplicitous politician could carefully engineer a gerrymander to satisfy the court's particular notion of compactness while still accomplishing the desired electoral outcome.  Throughout this article we will see many manifestations of this tension between the mathematical aim of precisely defining terms and the legal aim of providing healthy flexibility through careful use of imprecision.

\section{A matter of metrics}

We begin where, arguably, geometry itself begins: the Pythagorean formula.  In 2005, New York's Court of Appeals, the highest court in the state, heard the following case.  An undercover cop caught the appellant, James Robbins, selling drugs near Holy Cross grade school in Manhattan.  At dispute wasn't whether Robbins was illegally selling drugs, but whether he was doing so within 1,000 feet of the school's property, since that would subject him to a much harsher punishment.  The defense had argued that the shortest walkable path from the sale location to the school is 764 feet along one street then 490 feet along another, so the total distance is 1,256 feet, safely outside the drug-free school zone.  The prosecution had countered that the straight-line distance as the crow flies, or in other words the length of the hypotenuse computed using the Pythagorean formula, is the relevant distance, and that is just over 907 feet, placing the sale within the drug-free zone.  During the appeals case, Robbins' lawyers colorfully countered this assertion by noting that ``crows do not sell drugs,'' but the Court of Appeals nonetheless unanimously ruled against him and upheld the original felony conviction.

At first glance this case appears to hinge entirely on the fact that the statute in question specifies a numerical value for the distance defining school zones but not the metric to which this distance refers: the prosecution saw the Euclidean metric ($\ell_2$ norm) as appropriate, whereas the defense wanted, aptly enough considering the location of the incident, the Manhattan metric ($\ell_1$ norm).  Mathematically, neither metric is more valid than the other and they both produce distinct but well-defined 1000-foot school zones.  But there is more to the story here.  The defense did not literally argue $\ell_1$ over $\ell_2$ (which is a mathematical argument admitting no clear winner), the defense argued for using a ``pedestrian-distance'' defined as the shortest length of a path between two points along which a person could safely and legally walk.  But the Court of Appeals wisely recognized a serious problem with this notion: it is dynamic.  Indeed, the Chief Judge's opinion asserts: 
\begin{quote}
Plainly, guilt under the statute cannot depend on whether a particular building in a person's path to a school happens to be open to the public or locked at the time of a drug sale. [...] Further, if the pedestrian method of measurement were used, a drug dealer could take deliberate steps to evade the reach of the statute by putting up obstructions to render the footpath to a nearby school beyond the 1000-foot limit.
\end{quote}

By considering the intent of the statue (``to circumscribe a fixed geographical area, without regard to whether that area might contain obstacles around which people might have to detour'') and related laws and cases, the court is able to rule out the pedestrian-distance as a reasonable choice of metric.  But then if pedestrian mobility is irrelevant, there is scant reason to choose the Manhattan metric over the Euclidean metric, particularly since it was never explicitly argued for.  Thus the ruling by the Court of Appeals is actually quite astute and logical, rather than a mere mathematical whim.

One might think all this could simply have been avoided if the drug-free school zone statute had explicitly defined the intended metric.  In mathematics, a vaguely worded proposition is inherently problematic and precision of terms is of the utmost importance.  But it is impossible for the lawmaker to foresee the endlessly varied situations in which a statue will be applied, so small amounts of ambiguity, rather than being the result of careless oversight, can actually help ensure that the true intent of the law is followed.  It is conceivable, for instance, that when defining school zones, different metrics are appropriate in different situations and that forcing a Euclidean interpretation in all cases would lead to undesirable outcomes.  That said, too much ambiguity is also a problem, since that allows room for deliberate abuse of the law.  In fact, this was one of the main objections the Court of Appeals had to the pedestrian-distance: ``At a minimum, requiring speculation about pedestrian routes would create uncertainty in a statute which was meant to establish clear lines of demarcation.''  Thus, ironically, had the legislature been tempted to eliminate ambiguity by explicitly defining the school zones using pedestrian distance, it would have unintentionally created more ambiguity instead.

\section{The pancake theorem}

The next case involves a significant increase in the level of geometric sophistication over the preceding Pythagorean formula case---though here the geometric insight did not occur in the legal arguments or court ruling, but in an academic paper analyzing the case after the fact.  We shall discuss the 2015 U.S. Supreme Court case \emph{Evenwel v. Abbott} and Edelman's insightful response to it \cite{Edelman}.  This case concerns districting, which is also the context of gerrymandering, but the central issue here only tangentially relates to gerrymandering.

The Supreme Court ruled in 1964 that state legislatures must ensure their congressional districts have roughly equal populations.  But the ambiguity here is what exactly the term ``population'' is taken to mean.  While the predominant interpretation is total population, the courts often speak in terms of voting power, so only counting voters in the population does not seem unreasonable.  In \emph{Evenwel}, the appellants wanted Texas to switch from total population to citizen voting-age population in order to (in their view) equalize voting power thereby adhering to the intentions of the ``one person, one vote'' principle.  This distinction is especially salient in Texas, where there are large numbers of undocumented immigrants who are unable to vote and therefore would lose political power if they were not counted in the districting process.  As Edelman points out, the debate in \emph{Evenwel}, and in many years of academic literature on the issue more generally, has focused on the relative merits of these two populations (total versus citizen voting-age)---implicitly assuming that it is not possible to create districts that equalize both construals of population.  However, Edelman demonstrates not only that it \emph{is} possible to equalize both populations, but a mechanism for achieving this is a well-known result in topology.\footnote{Not just well-known: Francis Su, former president of the MAA, declared it one of his favorite theorems from topology! \cite{Su}}

The Ham Sandwich Theorem states that given $n$ bounded measurable subsets $S_1,S_2,\ldots,S_n \subset \mathbb{R}^n$ of $n$-dimensional Euclidean space, there exists an $(n-1)$-dimensional hyperplane that divides each $S_i$ into two subsets of equal measure.  In other words, we can cut $n$ blobs in half with a single linear slice.  There is no assumption that the subsets are disjoint, they may overlap or even be nested.  The name comes from the $n=3$ case where one pictures $A_1$ and $A_3$ as pieces of bread and $A_2$ as a piece of ham: a single slice cuts these all in half no matter how they are positioned (that is, even if they are not yet neatly assembled into a sandwich).  The $n=2$ case is called the Pancake Theorem, since one pictures two pancakes on a plate; this, as we shall see, is the relevant case for \emph{Evenwel}.  

The basic idea is to let the first pancake, $A_1$, be the set of all people in Texas, and let the second pancake, $A_2$, be the set of voting-age citizens in Texas (so $A_2 \subset A_1$).  Technically, since these are discrete sets we should either smooth each person out so that they have well-defined areas, or we could use the counting measure; since these result in essentially the same outcomes, let us sweep this detail under the rug. (The only real issues with counting discrete people is that an odd number clearly cannot be divided into two equal integers, but a $\pm 1$ error is insignificant in districting, and a dividing hyperplane may slice through a person, but that can be avoided by perturbing the locations very slightly.)  The Pancake Theorem then tells us that a single line through Texas cuts both populations evenly in half.  

We are done if we only need two districts, but how do we get more?  Edelman tells us.  By applying the same argument to each of our two districts (instead of the entire state of Texas), we can find two line segments that slice our two districts into four so that all four have equal total populations and equal citizen voting-age populations.  Repeating this construction inductively, we can equalize these two populations over $n=2^i$ polyhedral districts for any $i \ge 1$.  But what if the desired number of districts $n$ is not a power of two?  Dyadic fractions (numbers of the form $\frac{a}{2^b}$ for $a,b\in\mathbb{Z}$) are dense in the reals, indeed any $x\in\mathbb{R}$ can be made arbitrarily close to the dyadic fraction $\frac{\lfloor 2^i x\rfloor}{2^i}$ by taking $i$ sufficiently large, so if we want $n$ districts then we proceed as follows: 
\begin{enumerate}
\item fix $i\in\mathbb{N}$ and use the above procedure to produce $2^i$ districts equalizing both populations;
\item define $n$ new districts by grouping these $2^i$ districts into collections of size $\lfloor \frac{2^i}{n}\rfloor$ or $\lceil \frac{2^i}{n}\rceil$.
\end{enumerate}
Since these $n$ districts in step (2) differ from each other by at most one of the smaller districts in step (1), and these smaller districts each have $\frac{1}{2^i}$ of each population, the deviations from equality among the populations of the $n$ districts can be made arbitrarily small by taking $i$ sufficiently large.  Of course, larger values of $i$ result in more complicated district boundaries, so it might be desirable to seek districting outcomes using as small a value of $i$ as possible.  In Edelman's analysis he focuses on the large $i$ limit, but constraining $i$ leads to some interesting mathematical problems at the interface of planar geometry and combinatorics, such as:
\begin{itemize}
\item For fixed values of $n$ and $i$, how many districting outcomes are possible from Edelman's procedure? (Note that even a single instance of the Pancake Theorem can have non-unique cutting lines, so there are typically many choices in producing the $2^i$ districts in step (1) and then of course many further choices when grouping them according to step (2).)
\item For a given upper bound on the population deviations, what is the minimal value of $i$ such that there exists an equal districting outcome from Edelman's procedure in which all $n$ districts are connected?  What about connected and simply connected? (This is a reasonable constraint from a political perspective, since districts are often required to be contiguous and sometimes required not to have holes.)
\end{itemize}
These questions depend on the distribution of both populations across the state, so one could further study how the answers to these questions vary as the two populations vary. 

In summary, then, the legal debates (both academic and in the courtroom) have been about which population should be equalized, but Edelman shows that any two populations (for instance, the two considered in \emph{Evenwel}) can be equalized simultaneously, and his method for doing this relies on classical topology.  Edelman's procedure involves many choices and understanding the space of possible districting outcomes resulting from his procedure leads to novel research problems in pure math that potentially have significant practical applications.

\section{A cartographic controversy}\label{sec:maps}

Our next topic is a political matter rather than a legal one, but it is close in spirit to the themes in the rest of this paper and the geometry involved is sufficiently intriguing that we could not resist including a brief section on it.  In March of 2017, the Boston Public School system began displaying a new map of the Earth in their classrooms.  Word of this development quickly spread beyond local news sources such as \emph{Boston Globe}, so that within days there were articles in \emph{The Atlantic}, \emph{NPR}, and even \emph{The Guardian} in the U.K. about this little change to Boston classrooms.  What could possibly be so controversial about a map?  The answer is a fascinating blend of colonialism and differential geometry.

Mathematically speaking, the aim of cartography is to construct a function \[S^2 \rightarrow [a,b]\times [c,d] \subset \mathbb{R}^2\] from the sphere to a rectangular subset of the plane that respects various geometric properties of interest.  Since the sphere has positive Gaussian curvature whereas the plane has curvature zero, such a function cannot be a local isometry; in other words, distances are necessarily distorted when making a map of the globe.  That said, one could ask for angles to be preserved (these are \emph{conformal} maps) or for areas to be preserved (these are \emph{equiareal} maps)---but we cannot have both simultaneously since a conformal equiareal map is automatically an isometry and hence preserves curvature.  So differential geometry places a serious and unavoidable constraint on the cartographer.

The most widely seen map (at least in the U.S.) is the Mercator projection, a conformal map dating back to 1569 that can be described informally as follows: inscribe the sphere in a cylinder and expand the sphere like a balloon until it fills an entire segment of the cylinder, then cut it and unroll it onto a plane.  This map sends meridians (circles of constant longitude) and parallels (circles of constant latitude) to parallel lines meeting at right angles.  This implies that a straight line segment on the Mercator map corresponds to the path a ship would take if it keeps the compass direction constant---so the Mercator map was designed for navigational purposes and its utility in this regard continues today.  However, area is distorted: the further an object is from the equator, the greater will be its size on the Mercator map.  For instance, Africa appears roughly the same size as Europe on a Mercator map when in fact it is nearly three times larger.  Perhaps you can now see where this is going.

The Gall-Peters projection, first created in 1855 by Gall but popularized in the 1970s by Peters, is a particular equiareal map.  This is the map the Boston schools now display in their classrooms (often in addition to the Mercator map) as part of an effort to ``decolonize the classrooms.''  While the navigational advantage of Mercator is undeniable, classroom maps are not used for navigation, they are used to show students what the world looks like---and the main objection to the Mercator map in this context is that it inflates the size of Colonial powers.  This controversy did not begin in Boston; it was stirred up by Peters in the 1970s and even featured as a plot point in an episode of the TV show \emph{The West Wing} that aired in 2000.  Mercator versus Gall-Peters---conformal versus equiareal---this may be the only time in history that Gaussian curvature and its differential-geometric consequences were the subject of international news and popular TV.  See \cite{Rice,maps} for more on the math and history of the concepts discussed in this section.

\section{Is law like Euclidean geometry?}

Scholars as early as Gottfried Leibniz (1646--1716) have advocated using a ``geometric paradigm" in legal reasoning, though the form and meaning of this has evolved over time \cite{Leibniz}. The approach of Western legal systems in the 17th and 18th centuries was primarily \emph{inductive}: each case was considered individually according to the facts and then general rules were extrapolated from the deliberations in the case.  Since prior rulings were not examined in this framework, many inconsistencies and discrepancies arose. In an effort to fix this problem, legal theorists eventually turned to the axiomatic \emph{deductive} method of Euclid and attempted to establish a model of law based on geometric proof.  But can this really be done?

Laurence Tribe is one of the most esteemed constitutional scholars in America.  He wrote \emph{the} textbook on constitutional law and was Barack Obama's research mentor at Harvard Law; perhaps a less well-known distinction is that he was a PhD student in pure mathematics at Harvard before turning to legal studies.  We mentioned earlier that Tribe was at the epicenter of a debate on mathematics in the courtroom that took place in the early 70s.  In a more recent scholarly work \cite{TribeBook}, coauthored with his law student at the time (and now law professor at Cornell) Michael Dorf, Tribe addresses the question stated above of whether an axiomatic system modeled on Euclidean geometry can be used as a basis for jurisprudence. Tribe and Dorf ultimately reach a negative answer---more accurately, they find that doing so does not free the jurist from subjective value judgements---but their reasoning here is quite informative.

Tribe and Dorf recall Lakatos' influential perspective that mathematics, rather than being an accumulation of proven truths, is a process by which proofs are made more rigorous as they are subjected to counterexamples and criticism.  While Lakatos illustrated his philosophy with Euler's formula $V-E+F=2$ for polyhedra, Tribe and Dorf aim to reach a broader audience by replacing this with the high school proposition that the sum of the angles of any triangle is $180^\circ$.  They provide a standard, elementary proof of this fact, and then suppose someone comes along with an alleged counterexample: a triangle drawn on the surface of a sphere whose angle sum is strictly greater than $180^\circ$.  The three approaches identified by Lakatos that mathematicians typically take to such a situation are: (1) monster-barring, (2) exception-barring, or  (3) lemma-incorporation.

In the first case, one simply denies that a triangle on a sphere is in fact a triangle (it is, instead, a monster!), so the result remains valid by an ad hoc and brute-force denial of challenges: anything that disobeys the conclusions of the proposition is defined to be an object outside the scope of the proposition statement.  In the second case, one accepts the counterexample and adjusts the proposition statement accordingly, though still in an ad hoc manner: for all triangles that are not on the surface of a sphere, the sum of the angles is $180^\circ$.  Of course, this method does not address the gap in the proof that led to the counterexample, so one should be concerned that there are other counterexamples lurking around the corner.  The third method is far superior: the existence of a counterexample implies that at least one step in the original proof is false, so find what additional properties are needed to fix the false step(s).  In the Tribe-Dorf triangle example, the faulty step is one that relies on Euclid's fifth postulate (which does not hold for the surface of a sphere), so the lemma-incorporator refines the proposition statement as follows: for all triangles on surfaces for which Euclid's fifth postulate holds, the sum of the angles is $180^\circ$.  Then if one is faced with a particular surface and wonders whether this triangle property holds, it suffices to prove a lemma stating that Euclid's fifth postulate holds for this surface.

Tribe and Dorf then proceed to translate Lakatos' three approaches to the legal realm, and they discuss the pitfalls of each.  They provide interesting examples of real cases illustrating these ideas, but for the sake of brevity we shall only outline them in the abstract.  Tribe and Dorf identify monster-barring with the legal principle of drawing a distinction without a difference.  This means, for instance, that courts sometimes wish to grant a claimant standing without setting a broad precedent that would obligate them to do so for many related suits.  Thus, the court may refuse to hear a case that seems to have similar grounds as one they previously heard, and the justification for refusing the latter case is simply to deny the evident similarity: yes, we heard a case concerning triangles, but your spherical triangle is no triangle at all and so need not be granted standing.  Since the next approach, exception-barring, entails adjusting the proposition statement without concern for the purported proof, Tribe and Dorf draw a legal analogue to courts relying on the holdings of prior cases while ignoring their rationales.  The long-term problems this approach causes in the legal profession are quite similar those it causes in the mathematical profession, so we shall not elaborate further on it.  

But what of lemma-incorporation (which is, after all, the healthiest approach to take in mathematics)?  Unfortunately, the analogy between math and law breaks down here, at least according to Tribe and Dorf.  Let us quote directly from their book:
\begin{quote}
[...] mathematics, by definition, proceeds from assumed unprovable postulates.   Modern mathematicians do not argue about whether Euclid's Fifth Postulate is true in some metaphysical sense.  They know that some conjectures will be provable if that postulate is assumed true, and others will be provable if it is not.  And that is about all there is to say in the realm of mathematics.  By contrast, legal arguments center around the truth or falsity of the preliminary assumptions. [...] Law is, ultimately, unlike mathematics. \cite[p. 96]{TribeBook}
\end{quote} 
The main point, which Tribe and Dorf aptly illustrate with a historical example, is that one cannot go back to the legal ``proof'' of a result and find the false step; instead, one may disagree morally with a step of judicial reasoning, but there is no way to strengthen the legal argument by requiring an additional condition, like we do in mathematics.  A mathematician can say that the triangle sum formula holds for all Euclidean surfaces and make no claims about non-Euclidean surfaces, but a court cannot analogously say simply that abortion is protected by the Constitution for all those who find the act morally acceptable.

Curiously, a relation between legal foundations and non-Euclidean geometry was also mentioned by the famous mathematician Lipman Bers (1914--1993) in a separate context \cite[p. 199]{Calculus}.  He gave a lecture in which he compared the 1776 Declaration of Independence with Lincoln's 1863 Gettysburg Address by contrasting the unequivocal language of the former, ``We hold these truths to be self-evident, that all men are created equal,'' with the qualified wording of the latter: ``[...] conceived in liberty and dedicated to the proposition that all men are created equal.''  Why had this equality been demoted from self-evident (that is, axiomatic) to a proposition (that is, dependent upon axioms)?  Bers playfully suggested that this is answered by the advent of non-Euclidean geometry in 1830, and hence the recognition that the truth of a statement depends upon the choice of a system of axioms.

\section{Thinking like a topologist}

In the preceding section we summarized a chapter from a book by the mathematician-turned-lawyer Laurence Tribe where it is argued that mathematics is unlike law, at least in regard to founding the discipline on axiomatic deduction.  But that does not imply that mathematical training, or even mathematical thinking, has no jurisprudential utility.  In personal correspondence with the first author, Tribe had this to say about how he parlayed his graduate studies in mathematics into his legal studies:
\begin{quote}
Making the adjustment to the less precise and more ambiguously framed
questions, issues, and arguments on which legal thought is based took me a
couple of years, but in the end I was very glad to have made the switch. As
a legal analyst and writer, both in my scholarly work (including my
teaching) and in my advocacy, I find that my math background suggests novel,
and often useful, perspectives on problems that others tackle less
effectively for want of such a background. It's the sensitivity to matters
of structure and not anything quantitative (more topology than calculus, as
it were) that I find entails the largest transfer from mathematical to legal
thought.
\end{quote}

Tribe is not the only scholar to find topology in legal thought.  Meyerson in the book \emph{Political Numeracy} \cite{Numeracy} notes that the federalist and constitutional structure of the U.S. determines the shape of our national governance, by balancing the power of the Union with that of the States, but that within this shape there is enormous flexibility that has allowed our government to develop and grow over the years.  Meyerson analogizes this with the mathematical structure of a manifold in which the geometry can be deformed drastically without altering the topological type.  To illustrate with a more concrete example, he describes the following:
\begin{quote}
Imagine the federal government as a sphere containing all of the powers granted by the Constitution.  Within that sphere, however, there exists a hole, consisting of the powers that are, in the words of the Tenth Amendment, ``reserved to the States.''  Over time, the size of the hole has grown and shrunk relative to the size of the sphere, but the hole must remain if the Constitution's topological structure is to remain intact.  
\cite[p. 138]{Numeracy}
\end{quote}

Continuing this tour of topological and geometric thinking in law---as opposed to the more precise manifestations of geometric concepts in the courtroom described earlier in this paper---we come now to the social scientist Donald Black.  In an interview titled ``The geometry of law,'' Black describes some instances where the idea of metric geometry, in a very vague and philosophical sense, can help illuminate power structures and inequalities in the courtroom.  He uses the idea of a multidimensional ``social space'' to situate legal conflicts:
\begin{quote}
Consider the geometry of a homicide: When someone kills a stranger, for example, the relational distance covered by the killing is greater than when someone kills a friend or relative. Every killing also has a vertical structure. If the killer is, say, an unemployed and impoverished member of the victim's family while the victim is the family's prosperous patriarch, the killing has an upward direction (from a lower to a higher social elevation) while the direction of the legal case is downward (against a defendant below the victim). These relational and vertical characteristics are just two of many elements that constitute the multidimensional structure of any case of conflict. Moreover, the multidimensional structure of each case predicts and explains how it will be handled --- how law will behave from one case to the next.
\cite{Black}
\end{quote}
While it is impossible to precisely define distance between the real-world concepts Black envisions, nor would it make sense to ask whether the triangle inequality holds, the deliberately vague use of geometric ideas employed here by Black does help one recognize subtle jurisprudential issues that are otherwise difficult to articulate.  In our opinion Black over-reaches when he claims with bizarre specificity that ``Law is a curvilinear function of relational distance,'' but that does not diminish the overall utility of his perspective.

\section{Conclusion}

Statistics has been used in the courtroom since the 19th century but, unsurprisingly, instances of geometry and topology in the legal realm are much more rare.  Nonetheless, in this paper we have seen a case where the defendant's punishment rests on the choice of a metric; a case where the court debated two supposedly competing goals without realizing the Pancake Theorem allows them both to be simultaneously satisfied; and a geo-political controversy stemming from the fact that conformal equiareal maps preserve Gaussian curvature.  We have also seen several instances where geometric reasoning does and does not apply to jurisprudence.  Underlying many of these topics is a tension between the mathematician's need for precision and the lawyer's need for the flexibility that stems from vagueness.  Another central issue is that while both math and law are predicated on logical reasoning, in math the axioms are clearly defined whereas in law they are implicit and subjective.  Regardless, the courtroom has provided an interesting setting in which a few abstract geometric notions have been discussed and debated and contemplated, and knowledge of these real-world events may well enrich the mathematical education of geometry students.

\bibliographystyle{alpha}
\bibliography{bib}

\end{document}